\documentclass[11pt]{article}
\usepackage[russian]{babel}
\usepackage[cp1251]{inputenc}
\usepackage{amsmath,latexsym}
\usepackage[psamsfonts]{amssymb}
\usepackage[dvips]{graphicx}
\usepackage[mathcal]{euscript}
\begin{document}
\begin{center}
\textbf{ Нахождения общего решения одного\\ вырождающегося уравнения с дробной производной }\\
\textbf{Б.Ю.Иргашев}\\
Наманганский инженерно-строительный институт, Узбекистан\\
Институт Математики им.В.И.Романовского АН РУз.\\
bahromirgasev@gmail.com
\end{center}
УДК.517.926.4

\textbf{Аннотация.} В статье построено общее решение одного уравнения с обобщенной производной Хилфера, имеющего вырождение. Частные решения представлены через функцию Килбаса-Сайго. Полученно представление решения начальной задачи типа задачи Коши.

 \textbf{Kлючевые слова.} Производная дробного порядка , вырождение, ряд, функция Килбаса-Сайго, задача Коши.

В последнее время специалистами интенсивно изучаются уравнения c участием производных дробного порядка с переменными коэффициентами. К числу таких уравнений относятся вырождающиеся уравнения.  В работе  [1] изучалось уравнение
\[D_{0x}^\alpha {t^\beta }u\left( t \right) = \lambda u\left( x \right),\,\,0 < x < b,\]	
где $0 < \alpha  < 1,\,\,\lambda  - $ спектральный параметр, $\beta  = const \ge 0.$
В работе [2] были найдены решения в замкнутой форме уравнений дробного порядка
\[\left( {D_{0 + }^\alpha y} \right)\left( x \right) = a{x^\beta }y\left( x \right) + f\left( x \right)\left( {0 < x < d \le \infty ,\alpha  > 0,\beta  \in R,a \ne 0} \right),\]
\[\left( {D_ - ^\alpha y} \right)\left( x \right) = a{x^\beta }y\left( x \right) + f\left( x \right)\left( {0 \le d < x < \infty ,\alpha  > 0,\beta  \in R,a \ne 0} \right),\]
с дробными производными Римана-Лиувилля на полуоси $\left( {0,\infty } \right)$ [3]
\[\left( {D_{0 + }^\alpha y} \right)\left( x \right) = {\left( {\frac{d}{{dx}}} \right)^{[\alpha ] + 1}}\frac{1}{{\Gamma \left( {1 - \{ \alpha \} } \right)}}\int\limits_0^x {\frac{{y\left( t \right)dt}}{{{{\left( {x - t} \right)}^{\{ \alpha \} }}}}} ,\left( {x > 0;\alpha  > 0} \right),\]
\[\left( {D_ - ^\alpha y} \right)\left( x \right) = {\left( { - \frac{d}{{dx}}} \right)^{[\alpha ] + 1}}\frac{1}{{\Gamma \left( {1 - \{ \alpha \} } \right)}}\int\limits_x^\infty  {\frac{{y\left( t \right)dt}}{{{{\left( {t - x} \right)}^{\{ \alpha \} }}}}} ,\left( {x > 0;\alpha  > 0} \right),\]
($[\alpha ]$  и $\{ \alpha \} $  означают целую и дробную части действительного числа $\alpha $ ). К таким уравнениям приводят прикладные задачи [4]. Пример такого уравнения дает уравнение теории полярографии [5]
\[\left( {D_{0 + }^{{1 \mathord{\left/
 {\vphantom {1 2}} \right.
 \kern-\nulldelimiterspace} 2}}y} \right)\left( x \right) = a{x^\beta }y\left( x \right) + {x^{ - {1 \mathord{\left/
 {\vphantom {1 2}} \right.
 \kern-\nulldelimiterspace} 2}}},\left( {0 < x, - {1 \mathord{\left/
 {\vphantom {1 {2 < \beta  \le }}} \right.
 \kern-\nulldelimiterspace} {2 < \beta  \le }}0} \right),\]
 возникающее при $a =  - 1$  в задачах диффузии [5].

В данной работе мы построим фундаментальную систему решений вырождающегося уравнения c обобщенной дробной производной Хилфера   и получим представление решения начальной задачи типа задачи Коши. Рассмотрим следующее уравнение
\[D_{0y}^{\left( {\alpha ,\beta } \right)\mu }u\left( y \right) = \lambda {y^m}u (y),y>0, \lambda  \in C, m \ge 0,\eqno(1)\]
где
\[m + \mu \left( {\alpha  - \beta } \right) \ge 0,\]
\[D_{0y}^{\left( {\alpha ,\beta } \right)\mu } = I_{0y}^{\mu \left( {i - \alpha } \right)}\frac{{{d^i}}}{{d{y^i}}}I_{0y}^{\left( {1 - \mu } \right)\left( {i - \beta } \right)} = D_{0y}^{ - \mu \left( {i - \alpha } \right)}D_{0y}^{i - \left( {1 - \mu } \right)\left( {i - \beta } \right)},\eqno(2)\]
\[I_{0y}^\alpha g\left( y \right) = \frac{1}{{\Gamma \left( \alpha  \right)}}\int\limits_0^y {\frac{{g\left( t \right)dt}}{{{{\left( {y - t} \right)}^{1 - \alpha }}}}} ,\]
\[i - 1 < \alpha ,\beta  < i,\,i = 1,2,...,0 \le \mu  \le 1.\]
Краевые и начальные задачи для уравнений с участием оператора вида (2) исследовались в работах [6-8], вырождающиеся уравнение с дробным производным Хилфера в работах [9-10], а с дробными производными Римана-Лиувилля и Капуто в работах [2,11-12].

 Решение уравнения (1) будем искать в виде
\[u\left( y \right) = \sum\limits_{k = 0}^\infty  {{c_k}{y^{ak + b}}},\eqno(3)\]
где ${c_k},a>0,b$ пока неизвестные вещественные числа.\\
Предварительно вычислим следующее:
\[D_{0y}^{i - \left( {1 - \mu } \right)\left( {i - \beta } \right)}{y^\delta } = \frac{1}{{\Gamma \left( {\left( {1 - \mu } \right)\left( {i - \beta } \right)} \right)}}{\left( {\frac{d}{{dy}}} \right)^i}\int\limits_0^y {\frac{{{t^\delta }}}{{{{\left( {y - t} \right)}^{1 - \left( {1 - \mu } \right)\left( {i - \beta } \right)}}}}} dt = \]
\[ = \frac{{\Gamma \left( {\delta  + 1} \right){{\overline {\left( {\delta  + \left( {1 - \mu } \right)\left( {i - \beta } \right)} \right)} }_i}}}{{\Gamma \left( {\delta  + 1 + \left( {1 - \mu } \right)\left( {i - \beta } \right)} \right)}}{y^{\delta  + \left( {1 - \mu } \right)\left( {i - \beta } \right) - i}},\]
где $\delta  >  - 1,\,{\overline {\left( a \right)} _i} = a\left( {a - 1} \right)...\left( {a - i + 1} \right).$\\
Далее
\[D_{0y}^{ - \mu \left( {i - \alpha } \right)}{y^{\delta  + \left( {1 - \mu } \right)\left( {i - \beta } \right) - i}} = \frac{1}{{\Gamma \left( {\mu \left( {i - \alpha } \right)} \right)}}\int\limits_0^y {\frac{{{t^{\delta  + \left( {1 - \mu } \right)\left( {i - \beta } \right) - i}}}}{{{{\left( {y - t} \right)}^{1 - \mu \left( {i - \alpha } \right)}}}}} dt = \]
\[ = \frac{{{t^{\delta  + \left( {1 - \mu } \right)\left( {i - \beta } \right) - i + \mu \left( {i - \alpha } \right)}}}}{{\Gamma \left( {\mu \left( {i - \alpha } \right)} \right)}}\int\limits_0^1 {\frac{{{t^{\delta  + \left( {1 - \mu } \right)\left( {i - \beta } \right) - i}}}}{{{{\left( {1 - t} \right)}^{1 - \mu \left( {i - \alpha } \right)}}}}} dt = \]
\[ = \frac{{\Gamma \left( {\delta  + \left( {1 - \mu } \right)\left( {i - \beta } \right) - i + 1} \right)}}{{\Gamma \left( {\delta  + \left( {1 - \mu } \right)\left( {i - \beta } \right) - i + 1 + \mu \left( {i - \alpha } \right)} \right)}}{t^{\delta  + \left( {1 - \mu } \right)\left( {i - \beta } \right) - i + \mu \left( {i - \alpha } \right)}},\]
при условии
\[\delta  + \left( {1 - \mu } \right)\left( {i - \beta } \right) - i >  - 1.\eqno(4)\]
Итак имеем формулу
\[D_{0y}^{\left( {\alpha ,\beta } \right)\mu }{y^\delta } = \frac{{\Gamma \left( {\delta  + 1} \right){{\overline {\left( {\delta  + \left( {1 - \mu } \right)\left( {i - \beta } \right)} \right)} }_i}}}{{\Gamma \left( {\delta  + 1 + \left( {1 - \mu } \right)\left( {i - \beta } \right)} \right)}} \times \]
\[\frac{{\Gamma \left( {\delta  + 1 + \left( {1 - \mu } \right)\left( {i - \beta } \right) - i} \right)}}{{\Gamma \left( {\delta  + 1 + \left( {1 - \mu } \right)\left( {i - \beta } \right) - i + \mu \left( {i - \alpha } \right)} \right)}}{y^{\delta  + \left( {1 - \mu } \right)\left( {i - \beta } \right) - i + \mu \left( {i - \alpha } \right)}}.\eqno(5)\]
Теперь подставим (3) в (1), затем используя формулу (5) получим формальное равенство
\[\sum\limits_{k = 0}^\infty  {{c_k}} \frac{{\Gamma \left( {ak + b + 1} \right){{\overline {\left( {ak + b + \left( {1 - \mu } \right)\left( {i - \beta } \right)} \right)} }_i}}}{{\Gamma \left( {ak + b + 1 + \left( {1 - \mu } \right)\left( {i - \beta } \right)} \right)}} \times \]
\[\frac{{\Gamma \left( {ak + b + 1 + \left( {1 - \mu } \right)\left( {i - \beta } \right) - i} \right)}}{{\Gamma \left( {ak + b + 1 + \left( {1 - \mu } \right)\left( {i - \beta } \right) - i + \mu \left( {i - \alpha } \right)} \right)}}{y^{ak + b + \left( {1 - \mu } \right)\left( {i - \beta } \right) - i + \mu \left( {i - \alpha } \right)}} = \]
\[ = \lambda \sum\limits_{k = 0}^\infty  {{c_k}} {y^{ak + b + m}}.\]
Пусть
\[a = m + \alpha \mu  + \beta \left( {1 - \mu } \right),\]
\[b = s - \left( {1 - \mu } \right)\left( {i - \beta } \right),s = 0,1,...,i - 1,\]
тогда используя равенство:
\[{\overline {\left( {ak + b + \left( {1 - \mu } \right)\left( {i - \beta } \right)} \right)} _i}\Gamma \left( {ak + b + 1 + \left( {1 - \mu } \right)\left( {i - \beta } \right) - i} \right) = \]
\[ = {\overline {\left( {ak + b + \left( {1 - \mu } \right)\left( {i - \beta } \right)} \right)} _{i - 1}}\Gamma \left( {ak + b + 2 + \left( {1 - \mu } \right)\left( {i - \beta } \right) - i} \right) = \]
\[ = \Gamma \left( {ak + b + 1 + \left( {1 - \mu } \right)\left( {i - \beta } \right)} \right),\]
получим
\[\sum\limits_{k = 0}^\infty  {{c_k}} \frac{{\Gamma \left( {ak + b + 1} \right){y^{\left( {k - 1} \right)a}}}}{{\Gamma \left( {ak + b + 1 + \left( {1 - \mu } \right)\left( {i - \beta } \right) + \mu \left( {i - \alpha } \right) - i} \right)}} = \lambda \sum\limits_{k = 0}^\infty  {{c_k}} {y^{ak}}.\]
Найдем неизвестные коэффициенты $c_{k}$
\[{c_k} = \frac{{\lambda {c_{k - 1}}\Gamma \left( {ak + b - \left( {\beta  + \mu \left( {\alpha  - \beta } \right)} \right) + 1} \right)}}{{\Gamma \left( {ak + b + 1} \right)}} = \]
\[= {\lambda ^k}{c_0}\prod\limits_{j = 0}^{k - 1} {\frac{{\Gamma \left( {a\left( {j + 1} \right) + b - \left( {\beta  + \mu \left( {\alpha  - \beta } \right)} \right) + 1} \right)}}{{\Gamma \left( {a\left( {j + 1} \right) + b + 1} \right)}}}.\]
Заметим,что
\[a\left( {j + 1} \right) + b - \left( {\beta  + \mu \left( {\alpha  - \beta } \right)} \right) + 1 \ge a + b - \beta  - \mu \left( {\alpha  - \beta } \right) + 1 \ge \]
\[ \ge m + \beta  + \mu \left( {\alpha  - \beta } \right) - \left( {1 - \mu } \right)\left( {i - \beta } \right) - \beta  - \mu \left( {\alpha  - \beta } \right) + 1 = \]
\[ = m - \left( {1 - \mu } \right)\left( {i - \beta } \right) + 1 > m - \left( {1 - \mu } \right) + 1 = m + \mu  \ge 0 \Rightarrow \]
\[a\left( {j + 1} \right) + b - \left( {\beta  + \mu \left( {\alpha  - \beta } \right)} \right) + 1 > 0.\]
Значит имеем следующее семейство линейно независимых решений
\[{u_s}\left( y \right) = {y^{s - \left( {1 - \mu } \right)\left( {i - \beta } \right)}}\sum\limits_{k = 0}^\infty  {{c_k}} {\left( {\lambda {y^{m + \alpha \mu  + \beta \left( {1 - \mu } \right)}}} \right)^k},\eqno(6)\]
где
\[s = 0,1,...,i - 1,\]
\[{c_k} = \prod\limits_{j = 0}^{k - 1} {\frac{{\Gamma \left( {\left( {\left( {\beta  + \mu \left( {\alpha  - \beta } \right)} \right)} \right)\left( {\frac{{aj}}{{\beta  + \mu \left( {\alpha  - \beta } \right)}} + \frac{{a + b}}{{\beta  + \mu \left( {\alpha  - \beta } \right)}} - 1} \right) + 1} \right)}}{{\Gamma \left( {\left( {\left( {\beta  + \mu \left( {\alpha  - \beta } \right)} \right)} \right)\left( {\frac{{aj}}{{\beta  + \mu \left( {\alpha  - \beta } \right)}} + \frac{{a + b}}{{\beta  + \mu \left( {\alpha  - \beta } \right)}}} \right) + 1} \right)}}} ,k = 1,2,...,\]
\[{c_0} = 1.\]
Представление (6) также можно записать в виде
\[{y^b}{{\rm{E}}_{\gamma ,\frac{a}{\gamma },\frac{{a + b}}{\gamma } - 1}}\left( {\lambda {y^a}} \right),\eqno(7)\]
где
\[\gamma  = \beta  + \mu \left( {\alpha  - \beta } \right) > 0,\]
\[{E_{\alpha ,m,l}}\left( z \right) = \sum\limits_{i = 0}^\infty  {{c_i}{z^i}} ,{c_0} = 1,{c_i} = \prod\limits_{j = 0}^{i - 1} {\frac{{\Gamma \left( {\alpha \left( {jm + l} \right) + 1} \right)}}{{\Gamma \left( {\alpha \left( {jm + l + 1} \right) + 1} \right)}}} ,i \ge 1\]
- функция Килбаса-Сайго (см.[2]).

Проверим выполнения условия (4). Т.к.
\[D_{0y}^{i - \left( {1 - \mu } \right)\left( {i - \beta } \right)}{y^b} = \frac{{\Gamma \left( {b + 1} \right){{\overline {\left( s \right)} }_i}}}{{\Gamma \left( {s + 1} \right)}}{y^{s - i}} = 0,\]
то нужно проверить условие
\[a + b + \left( {1 - \mu } \right)\left( {i - \beta } \right) - i >  - 1.\]
Имеем
\[a + b + \left( {1 - \mu } \right)\left( {i - \beta } \right) - i = \]
\[ = m + \mu \left( {\alpha  - \beta } \right) + s - \left( {i - \beta } \right) > \]
\[ > m + \mu \left( {\alpha  - \beta } \right) + s - 1 \ge \]
\[ \ge m + \mu \left( {\alpha  - \beta } \right) - 1 \ge  - 1,\]
т.е. условие (4) выполняется.

Отметим, в случае $\alpha  = \beta ,$  имеем дробную производную в смысле Хилфера и семейство решений (7) преобразуются к виду
\[{u_s}\left( y \right) = {y^{s - \left( {1 - \mu } \right)\left( {i - \alpha } \right)}}{E_{\alpha ,\frac{m}{\alpha } + 1, \frac{{m + s - \left( {1 - \mu } \right)\left( {i - \alpha } \right)}}{\alpha }}}\left( {\lambda {y^{m + \alpha }}} \right),s = 0,1,...,i - 1.\]

Теперь рассмотрим следующую начальную задачу типа задачи Коши
\[\left\{ \begin{array}{l}
D_{0y}^{\left( {\alpha ,\beta } \right)\mu }u\left( y \right) = \lambda {y^m}u\left( y \right),y > 0,\lambda  \in C,\\
i - 1 < \alpha ,\beta  < i,0 \le \mu  \le 1,m \ge 0,m + \mu \left( {\alpha  - \beta } \right) \ge 0,\\
\mathop {\lim }\limits_{y \to  + 0} \frac{{{d^j}}}{{d{y^j}}}\left( {{y^{ - \left( {1 - \mu } \right)\left( {i - \beta } \right)}}u\left( y \right)} \right) = {\varphi _j},j = 0,1,...,i - 1.
\end{array} \right.\eqno(8)\]
Учитывая (6), (7) получим решение задачи (8) в виде
\[\begin{array}{l}
\sum\limits_{k = 0}^{i - 1} {\frac{{{\varphi _k}}}{{k!}}} {y^b}{{\rm{E}}_{\gamma ,\frac{a}{\gamma },\frac{{a + b}}{\gamma } - 1}}\left( {\lambda {y^a}} \right),\\
b = k - \left( {1 - \mu } \right)\left( {i - \beta } \right),\\
a = m + \beta  + \mu \left( {\alpha  - \beta } \right),\\
\gamma  = \beta  + \mu \left( {\alpha  - \beta } \right).
\end{array}\]

\begin{center}
Литература
\end{center}
1. Нахушев А.М. Дробное исчисление и его применение. - М.: Физматлит. 2003. - 272 c.\\
2. Килбас А. А., Сайго М.  Решение в замкнутой форме одного класса линейных дифференциальных уравнений дробного порядка. Дифференц. уравнения, 33 (2), 1997. c. 195 - 204.\\
3. Самко С. Г., Килбас А. А., Маричев О. И. Интегралы и производные дробного порядка и некоторые их приложения. Минск. Наука и техника. 1987. - 688 с.\\
4. Oldham К. В., Spanier J. The fractional calculus. New York; London. 1974.\\
5. Wiener K.  Wiss. Z. Univ. Halle Math. Natur. Wiss. R.  1983. 32 (1), 1983. pp. 41 - 46.\\
6. Bulavatsky V.M.. Closed form of the solutions of some boundary problems for anomalous diffusion equation with Hilfer’s generalized derivative. Cybernetics and Systems Analysis. 30(4), 2014, pp. 570-577\\
7.Karimov E., Ruzhansky M., Toshtemirov B.  Solvability of the boundary-value problem for a mixed equation involving hyper-Bessel fractional differential operator and bi-ordinal Hilfer fractional derivative. Mathematical Methods in the Applied Sciences. 41(1), 2023, pp. 54-77.\\
8. Karimov E. T., Toshtemirov B. H. Non-local boundary value problem for a mixed-type
equation involving the bi-ordinal Hilfer fractional differential operators. Uzbek Mathematical Journal 2021, Volume 65, Issue 2, pp.61-77. DOI: 10.29229/uzmj.2021-2-5\\
9. Restrepo, J. E., Suragan, D. (2021). Hilfer-type fractional differential equations with variable coefficients. Chaos, Solitons and Fractals, 150, 111146. doi:10.1016/j.chaos.2021.111146\\
10.Yuldashev T.K.,Kadirkulov B.J., Bandaliyev R.A. On a Mixed Problem for Hilfer Type Fractional Differential Equation with Degeneration.
Lobachevskii Journal of Mathematicsthis link is disabled, 2022, 43(1), pp. 263–274\\
10. B.Kh. Turmetov, B.J. Kadirkulov. On a problem for nonlocal mixed-type fractional order equation with degeneration,
Chaos, Solitons and Fractals,Volume 146, 2021, 110835, ISSN 0960-0779, https://doi.org/10.1016/j.chaos.2021.110835.\\
11. Smadiyeva A.G.   Well-posedness of the initial-boundary value problems for the
time-fractional degenerate diffusion equations. Bulletin of the Karaganda University. Mathematics series. 107(3), 2022, pp. 145-151.\\
12. Kilbas, Anatoly A.; Srivastava, Hari M.; Trujillo, Juan J. Theory and applications of fractional differential equations. North-Holland Mathematics Studies, 204. Elsevier Science B.V., Amsterdam, 2006. xvi+523 pp. ISBN: 978-0-444-51832-3; 0-444-51832-0

\end{document}